\documentclass[a4paper, reqno]{amsart}
\usepackage{amsmath, amssymb, eucal, amscd, amstext, enumerate}
\usepackage{amsfonts,amsthm}
\theoremstyle{plain}
\newtheorem{thm}{Theorem}[section]
\newtheorem{lem}[thm]{Lemma}
\newtheorem{eg}[thm]{Example}
\newtheorem{prop}[thm]{Proposition}
\newtheorem{cor}[thm]{Corollary}
\newtheorem{defn}[thm]{Definition}

\newtheorem{rem-ntn}[thm]{Remark and Notation}

\newenvironment{prf}{{\noindent \textbf{Proof:}\ }}{\hfill $\Box$\\ \smallskip}
\newenvironment{ack}{{\noindent \textbf{Acknowledgements} }}{\hfill \medskip}

\numberwithin{equation}{section}

\newcommand{\smnoind}{{\smallskip\noindent}}

\newcommand{\ti}{\tilde}

\newcommand{\id}{{\rm id}}

\newcommand{\norm}[1]{\left\|#1\right\|}
\newcommand{\sph}{\mathfrak{S}_1}
\newcommand{\botimes}{\bar{\otimes}}

\newcommand{\CL}{\mathcal{L}}

\newcommand{\KH}{\mathfrak{H}}
\newcommand{\KK}{\mathfrak{K}}
\newcommand{\KL}{\mathfrak{L}}
\newcommand{\KI}{\mathfrak{I}}

\newcommand{\CK}{\mathcal{K}}

\newcommand{\CA}{A}

\newcommand{\CHS}{\mathcal{HS}}

\newcommand{\BC}{\mathbb{C}}

\newcommand{\BN}{\mathbb{N}}

\newcommand{\BG}{{\mathbb{G}}}
\newcommand{\BH}{{\mathbb{H}}}

\newcommand{\Bo}{\mathbf{1}}
\newcommand{\BGd}{{\widehat{\mathbb{G}}}}
\newcommand{\Rep}{\mathrm{Rep}}
\newcommand{\ru}{\mathrm{u}}

\newcommand{\circledT}{\raisebox{.6pt}{\textcircled{\raisebox{-.5pt} {$\scriptstyle\top$}}}}

\newcommand{\C}[1]{\mathcal{#1}}

\begin{document}

\title{Property $T$ for general locally compact quantum groups}



\thanks{The authors are supported by the National Natural Science Foundation of China (11071126 and 11471168).}

\author{Xiao Chen \and Chi-Keung Ng}

\address[Xiao Chen]{Chern Institute of Mathematics, Nankai University, Tianjin 300071, China.}
\email{cxwhsdu@126.com}
\address[Chi-Keung Ng]{Chern Institute of Mathematics and LPMC, Nankai University, Tianjin 300071, China.}
\email{ckng@nankai.edu.cn; ckngmath@hotmail.com}

\begin{abstract}

In this short article, we obtained some equivalent formulations of property $T$ for a general locally compact quantum group $\BG$, in terms of
the full quantum group $C^*$-algebras $C_0^\ru(\widehat{\BG})$ and the $*$-representation of $C_0^\ru(\widehat{\BG})$ associated with the
trivial unitary corepresentation (that generalize the corresponding results for locally compact groups).
Moreover, if $\BG$ is of Kac type, we show that $\BG$ has property $T$
if and only if every finite dimensional irreducible $*$-representation of $C_0^\ru(\widehat{\BG})$ is an isolated point in the spectrum of
$C_0^\ru(\widehat{\BG})$ (this also generalizes the corresponding locally compact group result).
In addition, we give a way to construct property $T$ discrete quantum groups using bicrossed products.
\end{abstract}

\maketitle

\section{Introduction}

The notion of property $T$ for locally compact groups was first introduced by Kazhdan in the 1960s (see \cite{Kah}), and this property was proved to be a very useful notion.
A locally compact group $G$ is said to have \emph{property T} if every unitary representation of $G$ having almost invariant unit vectors actually has a non-zero invariant vector (see \cite[\S 1.1]{BHV}).
There are several equivalent formulations for property $T$ (see \cite{Val84} as well as Sections 1.1 and 1.2 of \cite{BHV}):
\begin{enumerate}[(P1)]
\item The full group $C^*$-algebra $C^*(G)$ can be decomposed as $\ker \pi_{\Bo_G} \oplus \BC$, where $\pi_{\Bo_G}$ is the $*$-representation induced by the trivial one dimensional representation $\Bo_G$.
\item There exists a minimal projection $p\in M(C^*(G))$ such that $\pi_{\Bo_G}(p) = 1$.
\item $\Bo_G$ is an isolated point in the topological space $\widehat{G}$ of irreducible unitary representations of $G$.
\item All finite dimensional elements in $\widehat{G}$ are isolated points in $\widehat{G}$.
\item There exists a finite dimensional element in $\widehat{G}$ which is an isolated point in $\widehat{G}$.
\end{enumerate}

\medskip

In \cite{Fim}, P. Fima extended the notion of property $T$ to discrete
quantum groups and he showed in \cite[Propositions 7 and 8]{Fim} that discrete quantum groups with property $T$  are of Kac type  and finitely generated (in some sense).
Recently, D. Kyed and M. Soltan studied property $T$ for discrete
quantum groups in \cite{KS}, using the techniques in the theory of matrix quantum groups, and they obtained the equivalences of property $T$ with the corresponding statements of (P1), (P2) and (P3) in the discrete quantum groups case.
They also proved that in the case when the discrete quantum group is unimodular (or equivalently, of Kac type), property $T$ is also equivalent to the corresponding statements of (P4) and (P5).
Furthermore, M. Daws, P. Fima, A. Skalski and S. White extended in \cite{DFSW} the definition of property $T$ to general locally compact quantum groups and showed that a locally compact quantum group has both the Haagerup property and property $T$ if and only if it is compact.

\medskip

Following the works of Fima, Kyed-Soltan as well as Daws-Fima-Skalski-White,  among others, the present article
devotes to the study of property $T$ for locally compact quantum groups.
We will extend the equivalences of property $T$ with (P1), (P2) and (P3) to locally compact quantum groups, and will verify that they are equivalent to the corresponding statements of (P4) and (P5) in the case of locally compact quantum groups of Kac type (which include all locally compact groups).
In fact, by employing $C^*$-algebras technique instead of matrix quantum groups technique (which do not work in this full generality), our proofs for these more general results are actually simpler than the ones in \cite{KS}.

\medskip

In Section 2, we will recall a basic fact on the spectra of $C^*$-algebras and will recall some notations and known facts on locally compact quantum groups.

\medskip

In Section 3, we will give a very short proof for the equivalences of property $T$ with the corresponding statements of (P1), (P2) and (P3) using the materials in Section 2.
On the other hand, in order to show the equivalence of property $T$ with the corresponding statements of (P4) and (P5), we need to consider contragredient unitary corepresentation of a unitay corepresentation.
We will use the technique in \cite{Ng}, concerning contragredient corepresentation of Kac algebras, to generalize \cite[Proposition A.1.12]{BHV} to the quantum case.
We then use it to get the desired equivalence.

\medskip

At the end of Section 3, we also present a new way to construct property $T$ discrete quantum groups.
Up to now, apart from the one given in \cite[Example 3.1]{Fim}, the only known examples of property $T$ discrete quantum groups are finite quantum groups and property $T$ discrete groups, as well as their direct products.
We will show how to construct property $T$ discrete quantum groups of Kac type using bicrossed products.

\bigskip

\section{Notations and preliminary}

\medskip

In this article, we use the convention that the inner product $\langle \cdot , \cdot \rangle$ of a complex Hilbert space $\KH$ is conjugate-linear in the first variable.
We denote by $\CL(\KH)$ and $\CK(\KH)$ the set of bounded linear operators and that of compact operators on $\KH$, respectively.
For any $x,y,z \in\KH$ and $T\in\CL(\KH)$, we denote by $\omega_{x,y}$ the normal functional given by $$\omega_{x,y}(T):=\langle x,Ty\rangle.$$
We set $\sph(\KH)$ to be the unit sphere of $\KH$.

\medskip

For a $C^*$-algebra $\CA$, we use $\Rep(\CA)$ to denote the collection of unitary equivalence classes of non-degenerate $*$-representations of $\CA$.
We consider $\widehat{\CA}\subseteq \Rep(\CA)$ to be the subset consisting of irreducible representations, equipped with the Fell topology (see \cite{Fell}). The topological space $\widehat \CA$ is known as the \emph{spectrum} of $A$.
Furthermore, all tensor products of $C^*$-algebras in the article, if not specified, are the minimal tensor products.

\medskip

Let us also recall some well-known facts concerning $\Rep(\CA)$ and $\widehat{\CA}$.
Suppose that $(\mu, \KH),(\nu, \KK)\in \Rep(A)$.
We write $\nu\subset \mu$ if there is an isometry $$V:\KK\to \KH$$ such that $$\nu(a) = V^*\mu(a)V \quad (a\in A).$$
Moreover, we write $\nu\prec \mu$ if $\ker \mu\subset \ker \nu$.

\medskip

The following lemma is well-known.
For the equivalence of Statements (1) and (2), one may use \cite[Theorems 1.2 and 1.6]{Fell}.
For the equivalence of Statements (1) and (3), one may use the fact that the topology on $\widehat{A}$ coincides with the one induced by the hull-kernel topology on ${\rm Prim}(A)$.
On the other hand, part (b) is a direct consequence of \cite[Lemma 1.11]{Fell}.

\medskip

\begin{lem}\label{lem32}
Let $A$ be a $C^*$-algebra and $(\mu,\KH)\in \widehat{A}$.

\smnoind
(a) The following statements are equivalent.
\begin{enumerate}
\item $\mu$ is an isolated point in $\widehat{A}$.


\item If $\pi\in \Rep(A)$ satisfying $\mu\prec \pi$, one has $\mu\subset \pi$.

\item $A = \ker \mu \oplus \bigcap_{\nu\in \widehat{A}\setminus \{\mu\}} \ker \nu$.
\end{enumerate}

\smallskip\noindent
(b) If $\dim \KH <\infty$, then $\{\mu\}$ is a closed subset of $\widehat{\CA}$.
\end{lem}

\medskip

Next, we recall some materials on locally compact quantum groups.
In the following, $(C_0(\BG), \Delta, \varphi, \psi)$ is a \emph{reduced $C^*$-algebraic locally compact quantum group} as introduced in
\cite[Definition 4.1]{KV1} (for simplicity, we will denote it by $\BG$).
The dual locally compact quantum group of $\BG$ (as defined in \cite[Definition 8.1]{KV1}) is denoted by $(C_0(\widehat{\BG}), \widehat{\Delta}, \widehat{\varphi}, \widehat{\psi})$.
We use $L^2(\BG)$ to denote the Hilbert space given by the GNS construction of the left invariant Haar weight $\varphi$ and consider both $C_0(\BG)$ and $C_0(\widehat{\BG})$ as $C^*$-subalgebras of $\CL(L^2(\BG))$.
There is a unitary $$W_\BG\in M(C_0(\BG)\otimes C_0(\BGd))\subseteq \CL(L^2(\BG)\otimes L^2(\BG)),$$ called the \emph{fundamental multiplicative unitary} that implements the comultiplication: $$\Delta(x)=W_\BG^*(1\otimes x)W_\BG \quad (x\in C_0(\BG)).$$
The von Neumann subalgebra $L^{\infty}(\BG)$ generated by $C_0(\BG)$
in $\CL(L^2(\BG))$ is a Hopf von Neumann algebra under a comultiplication $\tilde{\Delta}$ defined by $W_\BG$ as in the above (see \cite{KV2} or \cite[Section 8.3.4]{Tim}).

\medskip

\begin{defn}
For any Hilbert space $\KH$, a unitary $U\in M(\CK(\KH)\otimes C_0(\BG))$ is called a \emph{unitary corepresentation} of $\BG$ on $\KH$ if
\begin{equation}\label{eqt:defn-corep}
(\id\otimes\Delta)(U)=U_{12}U_{13},
\end{equation}
where $U_{ij}$ is the usual ``leg notation'' (see e.g.\ \cite{BS}).
\end{defn}

\medskip

The universal version of $\BGd$ is denoted by $(C_0^\ru(\BGd),\widehat{\Delta}^\ru)$ (see \cite[Section 4 and 5]{Kus}).
As shown in \cite[Proposition 5.2]{Kus}, there exists a unitary $$V_\BG^\ru\in M(C_0^\ru(\widehat{\BG})\otimes C_0(\BG))$$ that implements a bijection between unitary corepresentations $U$ of $\BG$ on $\KH$ and non-degenerate $*$-representations $\pi_U$ of $C_0^\ru(\widehat{\BG})$ on $\KH$ through the correspondence
$$U = (\pi_U\otimes\id)(V_\BG^\ru).$$
The identity $\Bo_\BG$ of $M(C_0(\BG))$ is a unitary corepresentation of $\BG$ on $\BC$ and $\pi_{\Bo_\BG}$ is a character of $C_0^\ru(\widehat{\BG})$.

\medskip

If $W$ is another unitary corepresentation of $\BG$ on a Hilbert space $\KK$, we denote by $U\circledT W$ the unitary corepresentation
$U_{13}W_{23}$ on $\KH\otimes\KK$ and call it the \emph{tensor product} of $U$ and $W$.
In this case,
\begin{equation}\label{eqt:ten-prod}
\pi_{U\circledT W} = (\pi_{U}\otimes\pi_{W})\circ \widehat{\Delta}^\ru.
\end{equation}

\medskip

\begin{defn}
Let $U\in M(\CK(\KH)\otimes C_0(\BG))$ be a unitary corepresentation.

\smnoind
(a) $\xi\in\KH$ is called a \emph{$U$-invariant vector} if for every $\eta\in L^2(\BG)$, one has
$$U(\xi\otimes\eta)=\xi\otimes\eta.$$

\smnoind
(b) A net $\{\xi_i\}_{i\in \KI}$ in the unit sphere $\sph(\KH)$ is called an \emph{almost $U$-invariant unit vector} if for each $\eta\in L^2(\BG)$, one has $$\norm{U(\xi_i\otimes\eta)-\xi_i\otimes\eta}\rightarrow 0.$$
\end{defn}

\medskip

The following proposition can be found in \cite[Theorem 5.1]{BT} and \cite[Proposition 2.7]{DFSW}.

\medskip

\begin{prop}\label{prop:inv-almost-inv}
Let $U$ be a unitary corepresentation $U$ of $\BG$ on $\KH$.

\smnoind
(a) An element $\xi \in \KH$ is $U$-invariant if and only if for all $x\in C^\ru_0(\widehat{\BG})$, one has
$$\pi_{U}(x)\xi=\pi_{\Bo_\BG}(x)\xi.$$

\smnoind
(b) A net $\{\xi_i\}_{i\in \KI}$ in $\sph(\KH)$ is an almost $U$-invariant unit vector if and only if for all $x\in C^\ru_0(\widehat{\BG})$, one has $$\norm{\pi_{U}(x)\xi_i-\pi_{\Bo_\BG}(x)\xi_i}\rightarrow 0.$$
\end{prop}

\medskip

As in the literature, we write $U\subset W$ and $U\prec W$ when $\pi_U\subset \pi_W$ and $\pi_U\prec \pi_W$, respectively (see, e.g., \cite[Section 5]{BT} or \cite[Definition 2.3]{DFSW}).
From Proposition \ref{prop:inv-almost-inv}, we can get directly the following corollary.

\medskip

\begin{cor}\label{cor26}
Let $U$ be a unitary corepresentation of $\BG$.

\smnoind
(a) $U$ has a non-zero invariant vector if and only if $\Bo_{\BG}\subset U$.

\smnoind
(b) $U$ has almost invariant vectors if and only if $\Bo_{\BG}\prec U$.
\end{cor}

\bigskip

\section{Property $T$ for locally compact quantum groups}

\medskip

\begin{defn}\label{defn31}
A locally compact quantum group $\BG$ is said to have \emph{property $T$} if every unitary corepresentation having an almost invariant unit vector has a non-zero invariant vector.
\end{defn}

\medskip

Let us first generalize the equivalences of property $T$ with the corresponding statements of (P1), (P2) and (P3) to the general case of locally compact quantum groups.
Note that our proof here is even simpler than the case of locally compact groups (by using the materials in Section 2).

\medskip

\begin{prop}\label{prop:first-equiv}
The following statements are equivalent.

\begin{enumerate}[(T1)]
\item $\BG$ has property $T$

\item $\pi_{\Bo_\BG}$ is an isolated point in $\widehat{C_0^\ru(\widehat{\BG})}$.

\item $C_0^\ru(\widehat{\BG}) = \ker \pi_{\Bo_\BG}\oplus \BC$.

\item There is a projection $p_{\BG} \in M(C_0^\ru(\widehat{\BG}))$ with $$p_{\BG}C_0^\ru(\widehat{\BG})p_{\BG} = \BC p_\BG\ \;\text{and}\ \; \pi_{\Bo_\BG}(p_{\BG})=1.$$
\end{enumerate}
\end{prop}
\begin{prf}
$(T1)\Leftrightarrow (T2)$
This follows from Corollary \ref{cor26} and Lemma \ref{lem32}(a).

\smnoind
$(T2)\Rightarrow (T3)$.
This follows from Lemma \ref{lem32}(a).

\smnoind
$(T3)\Rightarrow (T4)$.
One may take $p_{\BG}=(0,1)\in \ker \pi_{\Bo_\BG}\oplus \BC$.

\smnoind
$(T4)\Rightarrow (T2)$.
Let $x$ be an element in $C_0^\ru(\widehat{\BG})$ such that $\pi_{\Bo_\BG}(x)=1$. From $$p_{\BG}=p_{\BG}x p_{\BG},$$ we know that $p_{\BG}$
actually belongs to $C_0^\ru(\widehat{\BG})$.
As $p_{\BG}C_0^\ru(\widehat{\BG})p_{\BG}$ is a hereditary $C^*$-subalgebra of $C_0^\ru(\widehat{\BG})$, its spectrum can be identified as an open subset of $\widehat{C_0^\ru(\widehat{\BG})}$.
In fact, this open subset is $\{\pi_{\Bo_\BG}\}$, by \cite[Lemma 1]{Val84}.
On the other hand, by Lemma \ref{lem32}(b), $\{\pi_{\Bo_\BG}\}$ is also a closed subset of $\widehat{C_0^\ru(\widehat{\BG})}$.
\end{prf}

\medskip

\begin{eg}\label{eg:non-T}
Let $G$ be a locally compact group.

\smnoind
(a) Suppose that $G_1$ and $G_2$ are closed subgroup of $G$ such that the canonical map $\varphi:G_1\times G_2\to G$ is a bijective homeomorphism into an open dense subset $\Omega$ of $G$ (and hence $G\setminus \Omega$ has measure zero).
We consider $\alpha$ and $\beta$ to be canonical continuous actions of $G_1$ and $G_2$ on $G_1\backslash G$ and $G/G_2$, respectively.
Then $G_1$ and $G_2$ is a matched pair of locally compact groups in the sense of  \cite[Definition 3.6.7]{Vaes01}.
By considering the trivial cocycles, one obtained from \cite[Theorem 3.4.13]{Vaes01} the locally compact quantum group $\BG$.
In fact, the fundamental unitary of $\BG$ is the unitary $V$ as given in \cite{BS98} and one has $C_0(\BG) = C_0(G/ G_2)\rtimes_{\beta, r} G_2$ and $C_0(\widehat{\BG}) = C_0(G_1\backslash G)\rtimes_{\alpha, r} G_1$.

Now, suppose that $\Omega = G$, both $G_1$ and $G_2$ are amenable with $G_1$ being non-compact.
By \cite[Theorem 6]{Ng-Kac-sys}, $V$ is amenable, or equivalently, $\BGd$ is coamenable.
If $\BG$ has property $T$, then \cite[Proposition 6.2]{DFSW} implies that $\BG$ is compact and hence $C_0(\BG)$ is unital.
This gives the contradiction that  $G_1\cong G/G_2$ is compact.

\smnoind
(b) Suppose that $\widehat{\BG}$ is the dual group of the locally compact quantum group $\BG$ corresponding to $G$.
Since $C_0^\ru(\BG) = C_0(G)$ and $\Bo_\BG$ corresponding to the evaluation at the identity $e$ of $G$, we know from Proposition \ref{prop:first-equiv} that $\widehat{\BG}$ has property $T$ if and only if $G$ is discrete.
This part can also be deduced from \cite[Proposition 6.2]{DFSW}.
\end{eg}

\medskip

We recall that $\BG$ is said to be \emph{of Kac type} if $L^{\infty}(\BG)$ is a Kac algebra (see e.g.\ \cite{ES}).
In this case, the antipode is bounded.
We want to extend the equivalences of property $T$ with the corresponding statements of (P4) and (P5) in the Kac type case.
Before that we need to generalize \cite[Proposition A.1.12]{BHV} to this case.
Let us set some more notations.

\medskip

\emph{From now on, $\BG$ is of Kac type, $U$ and $V$ are unitary corepresentations of $\BG$ on $\KH$ and $\KK$ respectively.}

\medskip

One may regard $U\in\CL(\KH)\botimes L^{\infty}(\BG)$ and $V\in\CL(\KK)\botimes L^{\infty}(\BG)$.
As in \cite{Ng}, we define the contragredient $\bar V$ of $V$ by
$$\bar{V}:=(\tau\otimes\kappa)(V),$$
where $\tau$ is the canonical anti-isomorphism from $\CL(\KK)$ to $\CL(\bar{\KK})$ (with $\bar \KK$ being the conjugate Hilbert space of $\KK$) and $\kappa$ is the bounded antipode on $L^\infty(\BG)$.
Then $\bar V$ is a unitary corepresentation of $\BG$ on $\bar \KK$ (see, e.g., \cite[Corollary A.6(d)]{BS} or \cite[Remark 2.2]{Ng}).

\medskip

There is a canonical bijective isometry $\Theta$ from $\KH\otimes\bar \KK$ to the Hilbert space $\CHS(\KK,\KH)$ of Hilbert-Schmidt operators given by
$$\Theta(x\otimes \bar y)(z) := x \langle y,z\rangle,$$ for any $x\in\KH\ \;\text{and}\ \;\bar{y}\in\bar \KK$.
We set
$$(U; V)_{\CHS}\ :=\ (\Theta\otimes\id)(U\circledT \bar V)(\Theta^*\otimes \id),$$
which is a unitary corepresentation of $\BG$ on $\CHS(\KK,\KH)$.

\medskip

\begin{lem}\label{lem:invar-HS}
$T\in \CHS(\KK,\KH)$ is $(U; V)_\CHS$-invariant if and only if $$U(T\otimes 1)V^* = T\otimes 1.$$
\end{lem}
\begin{prf}
There are sequences $\{\xi_k\}_{k\in \BN}$ and $\{\eta_k\}_{k\in \BN}$ in $\KH$ and $\KK$, respectively, with
$$\Theta\left(\sum_{k\in \BN} \xi_k\otimes \bar \eta_k\right) = T.$$
By \cite[Lemma 3.8(a)]{Ng}, for any $\xi\in \KH$, $\eta\in \KK$ and $\alpha,\beta\in L^2(\BG)$, one has
\begin{eqnarray*}
\big\langle \xi\otimes \bar \eta \otimes \beta , (U\circledT \bar V)\big(\Theta^*(T)\otimes \alpha\big) \big\rangle
& = & {\sum}_{k\in \BN} \big\langle U_{13}^*(\xi\otimes \eta_k\otimes \beta) , V_{23}^*(\xi_k\otimes \eta\otimes \alpha)\big\rangle \\
& = & {\sum}_{k\in \BN} \big\langle U^*(\xi\otimes \beta) , \xi_k\otimes (\omega_{\eta_k,\eta} \otimes \id)(V^*) \alpha\big\rangle\\
& = & {\sum}_{k\in \BN} \big\langle \xi\otimes \beta , U\big(\Theta(\xi_k\otimes \bar \eta_k)\otimes 1\big)V^*(\eta \otimes \alpha)\big\rangle\\
&  = & \big\langle \xi\otimes \beta , U(T\otimes 1)V^*(\eta \otimes \alpha)\big\rangle,
\end{eqnarray*}
because  $\Theta(\xi_k\otimes \bar \eta_k)S\eta = \omega_{\eta_k,\eta}(S)\xi_k$  ($S\in \CL(\KK)$)
and
$$\big\langle \xi\otimes \bar \eta \otimes \beta , \Theta^*(T)\otimes \alpha \big\rangle = \big\langle \xi\otimes \beta , (T\otimes 1)(\eta \otimes \alpha)\big\rangle.$$
Thus, we have
$$(U\circledT \bar V)\big(\Theta^*(T)\otimes \alpha\big) = \Theta^*(T)\otimes \alpha \qquad (\alpha\in L^2(\BG))$$
if and only if
$U(T\otimes 1)V^* = T\otimes 1$.
\end{prf}

\medskip

\begin{prop}\label{prop:fin-dim-subrep}
Let $U$ and $V$ be as in the above.
Then $\Bo_{\BG}\subset U\circledT \bar V$ if and only if there is a finite dimensional unitary corepresentation $W$ such that $W\subset U$ and $W\subset V$.
\end{prop}
\begin{prf}
$\Rightarrow)$.
By Corollary \ref{cor26} and Lemma \ref{lem:invar-HS}, there is $T\in \CHS(\KK,\KH)\setminus \{0\}$ such that $U(T\otimes 1)V^* = T\otimes 1$.
The proof now preceeds as that of \cite[Proposition A.1.12]{BHV}.

More precisely, since $TT^*\in\CK(\KH)_+$, there exists $\lambda\in\sigma(TT^*)\setminus\{0\}$ with the corresponding eigenspace
$\C{E}_{\lambda}$ being finite dimensional.
It follows from $U(TT^*\otimes 1)U^*=TT^*\otimes 1$ that
$$TT^* \pi_U = \pi_U TT^*$$ and $\C{E}_{\lambda}$ is $\pi_U$-invariant.
Moreover, from the equalities
$$\norm{T^*\xi}^2
=\langle\xi,TT^*\xi\rangle
=\lambda\norm{\xi}^2
\quad (\xi\in\C{E}_{\lambda}),$$
we know that $\lambda^{\frac{1}{2}}T^*|_{\C{E}_{\lambda}}: \C{E}_{\lambda}\rightarrow T^*(\C{E}_{\lambda})$ is a bijective isometry.
Furthermore, as $$V(T^*\otimes 1)= (T^*\otimes 1)U$$ and $\C{E}_{\lambda}$ is $\pi_U$-invariant, we know that $T^*(\C{E}_{\lambda})$ is $\pi_{V}$-invariant
and $$\pi_U|_{\C{E}_{\lambda}}\cong\pi_{V}|_{T^*(\C{E}_{\lambda})}$$ under $\lambda^{\frac{1}{2}}T^*|_{\C{E}_{\lambda}}$.
Consequently, $W=(\pi_U|_{\C{E}_{\lambda}}\otimes\id)(V_\BG^\ru)$ is the finite dimensional corepresentation that is demanded.

\smnoind
$\Leftarrow)$.
Let $\KL$ be a finite dimensional Hilbert space and $W\in M(\CK(\KL)\otimes C_0(\BG))$.
By Lemma \ref{lem:invar-HS}, the identity operator $1\in \CHS(\KL,\KL)$ is
$(W; W)_\CHS$-invariant.
Thus, Corollary \ref{cor26} gives $$\Bo_{\BG}\subset W\circledT \bar W \subset U \circledT \bar V$$ as required.
\end{prf}

\medskip

The proof of the following theorem now follows from similar lines of argument as that of \cite[Theorem 1.2.5]{BHV}.
For completeness, we present the argument here.

\medskip

\begin{thm}\label{thm38}
Let $\BG$ be a locally compact quantum group of Kac type.
Then property $T$ of $\BG$ is also equivalent to the following statements.

\begin{enumerate}[(T1)]
\setcounter{enumi}{4}

\item Every finite dimensional irreducible representation of $C_0^\ru(\widehat{\BG})$ is an isolated point in  $\widehat{C_0^\ru(\widehat{\BG})}$.

\item $C_0^\ru(\widehat{\BG}) \cong B \oplus M_n(\BC)$ for a $C^*$-algebra $B$ and an $n\in \BN$.
\end{enumerate}
\end{thm}
\begin{prf}
$(T2)\Rightarrow (T5)$.
Let $\mu\in \widehat{C_0^\ru(\widehat{\BG})}$ and $\pi\in \Rep(C_0^\ru(\widehat{\BG}))$ such that $\mu$ is
finite dimensional and $\mu\prec\pi$.
If we set $$U := (\pi\otimes \id)(V_\BG^\ru)\quad \text{and} \quad V:= (\mu\otimes \id)(V_\BG^\ru),$$
then $V\prec U$.
Therefore, by \eqref{eqt:ten-prod} and Proposition \ref{prop:fin-dim-subrep}, one has
$$\Bo_{\BG}\subset V \circledT \bar V\prec U\circledT \bar V.$$
Hence, $\Bo_{\BG}\subset U\circledT \bar V$ by Lemma \ref{lem32}(a) and Proposition \ref{prop:first-equiv}.
This gives a unitary corepresentation $W$ with $W\subset U$ and $W\subset V$ (again, by Proposition \ref{prop:fin-dim-subrep}) and the irreducibility implies $V = W\subset U$ (or equivalently, $\mu\subset \pi$).
Now, Lemma \ref{lem32}(a) gives the required conclusion.

\smnoind
$(T5)\Rightarrow (T6)$.
This follows from Lemma \ref{lem32}(a).

\smnoind
$(T6)\Rightarrow(T2)$.
Let $\mu$ be the irreducible $*$-representation of $C_0^\ru(\widehat{\BG})$ corresponding to the summand $M_n(\BC)$ and denote
$$U := (\mu\otimes \id)(V_\BG^\ru).$$
As $\mu$ is finite dimensional, one has
$$(\mu\otimes \bar \mu)\circ \widehat{\Delta}^\ru = \mu_0\oplus \cdots \oplus \mu_n$$ for some $\mu_0,...,\mu_n\in \widehat{C_0^\ru(\BGd)}$.
By Proposition \ref{prop:fin-dim-subrep}, we may assume that $\mu_0 = \pi_{\Bo_\BG}$.
Moreover, Lemma \ref{lem32}(b) tells us that all such
$\{\mu_k\}$ are closed subset of $\widehat{C_0^\ru(\BGd)}$.
Suppose on the contrary that $\{\pi_{\Bo_\BG}\}$ is not open in $\widehat{C_0^\ru(\BGd)}$.
Then there is a net $\{\sigma_i\}_{i\in \KI}$ in $\widehat{C_0^\ru(\BGd)}\setminus \{\mu_0,...,\mu_n\}$ that converges to  $\pi_{\Bo_\BG}$.
Thus, $\pi_{\Bo_\BG}\prec\bigoplus_{i\in \KI}\sigma_i$, which implies $$\mu \prec \bigoplus_{i\in \KI} (\sigma_i\otimes \mu)\circ \widehat{\Delta}^\ru.$$
By Lemma \ref{lem32}(a), one has
$\mu \subset \bigoplus_{i\in \KI}(\sigma_i\otimes \mu)\circ \widehat{\Delta}^\ru$, and hence,
$$\mu \subset (\sigma_{i_0}\otimes \mu)\circ \widehat{\Delta}^\ru,$$ for some $i_0\in \KI$ (as $\mu$ is irreducible).
If we put
$$V_0 := (\sigma_{i_0}\otimes \id)(V_\BG^\ru),$$
then by Proposition \ref{prop:fin-dim-subrep}, we know that
$$\Bo_\BG \subset U\circledT \bar U \subset V_0\circledT U\circledT \bar U = {\bigoplus}_{k=0}^n V_0 \circledT (\mu_k\otimes \id)(V_\BG^\ru).$$
Now, Proposition \ref{prop:fin-dim-subrep} and the irreducibility of $\mu_k$ ($k=0,...,n$) again tells us that there is $k_0\in \{0,....,n\}$ with $$(\mu_{k_0}\otimes \id)(V_\BG^\ru) \subset V_0.$$
However, this will produce the contradiction that $\sigma_{i_0} = \mu_{k_0}$.
\end{prf}

\medskip

The following corollary is a direct consequence of Proposition \ref{prop:first-equiv} and Theorem \ref{thm38}.

\medskip

\begin{cor}\label{cor:quot}
	Let $\BH$ be another locally compact quantum group such that there is a surjective $^*$-homomorphism $\Phi:C_0^\ru(\widehat{\BG}) \to C_0^\ru(\widehat{\BH})$.
	Suppose that $\BG$ has property $T$.
	If either $\pi_{\Bo_\BG} = \pi_{\Bo_\BH}\circ \Phi$ or $\BG$ is of Kac type, then $\BH$ has property $T$.
\end{cor}

\medskip

Let us end this paper with a ``non-trivial'' construction of discrete quantum groups with property $T$.
Suppose that $G_1$ is a property $T$ discrete group acting non-trivially on a finite group $G_2$ by group automorphisms  and take $G$ to be the semi-direct product $G_2\rtimes G_1$.
For example, if $[G_1,G_1]\neq G_1$, then there exist a non-trivial action of the finite abelian group $G_1/[G_1,G_1]$ on some finite group $G_2$.
If $\beta$ is the action of $G_2$ on $G_1$ as in Example \ref{eg:non-T}(a), then $\beta$ is trivial.
Thus, the resulting locally compact quantum group $\BG$ in Example \ref{eg:non-T}(a) is of Kac type (see \cite[Corollary 3.6.17]{Vaes01}).

\medskip

Furthermore, the following theorem tells us that  $\BG$  has property $T$.
Observe that an essential part of the proof of this theorem is to verify that $C_0^u(\BG)$ is a quotient $C^*$-algebra of the full crossed product $C(G_2)\rtimes_\alpha G_1$.
This could be a known fact, but since we do not find it in the literature, we give an argument here for the sake of completeness.

\medskip

\begin{thm}
If $G$, $G_1$ and $G_2$ are as in the above, then the discrete quantum group $\BG$ as in Example \ref{eg:non-T}(a) has property $T$.
\end{thm}
\begin{prf}
Let $\alpha$ and $\beta$ be the actions as in Example \ref{eg:non-T}(a).
We denote by $\Delta$, $\Delta_1$ and $\hat \Delta_2$ the coproducts on $C_0(\BG)$, $C_0(G_1)$ and $C^*(G_2)$ respectively.
By abuse of notations, we use
$$\alpha: C_0(G_2)\to C_b(G_1\times G_2) \quad \text{and} \quad \beta:C_0(G_1)\to C_b(G_1\times G_2)$$
to denote the maps induced by the actions $\alpha$ and $\beta$ as in \cite[p.275]{Vaes01}, respectively.
In particular,
\begin{equation}\label{eqt:defn-alpha}
\alpha(y)(g,s)\ =\ y(\alpha_g(s)) \qquad (y\in C_0(G_2); g\in G_1;s\in G_2).
\end{equation}
The triviality of $\beta$ implies that $\beta(\phi) = \phi\otimes 1$ ($\phi\in C_0(G_1)$).
Suppose that $W^{(1)}\in M(C_0(G_1)\otimes C^*_r(G_1))$ and $\hat W^{(2)}\in M(C^*_r(G_2)\otimes C_0(G_2))$ are the fundamental unitary corresponding to $G_1$ and the dual of $G_2$, respectively.
As in \cite[Definition 3.4.2]{Vaes01}, the fundamental unitary of $\BG$ is given by
$$W_\BG = W^{(1)}_{13}\big(1 \otimes (\id\otimes \alpha)(\hat W^{(2)})\big)\in \CL\big(L^2(G_1\times G_2\times G_1\times G_2)\big).$$

Consider $\Bo_{G_2}$ to be the trivial reprensetation of $G_2$
 and $\eta_{e}: C_0(G_1)\to \BC$ to be the evaluation at the identity $e$ of $G_1$.
Since $G_2$ is finite and $\beta$ is trivial, the $C^*$-algebra $C_0(\BG)$ coincides with $C_0(G_1)\otimes C^*(G_2)$.
For any $z\in C_0(\BG)$, we have
\begin{eqnarray}\label{eqt:eta-coalg-homo}
\lefteqn{(\eta_e\otimes \id \otimes \eta_e\otimes \id)\Delta(z)}\quad \quad \quad \nonumber \\
& = & (\eta_e\otimes \id \otimes \eta_e\otimes \id)\big(W_\BG^*(1\otimes z)W_\BG\big)\nonumber \\
& = & \big(\id \otimes (\eta_e\otimes\id)\alpha\big)(\hat W^{(2)})^*\big(1\otimes (\eta_e\otimes \id)(z)\big)\big(\id \otimes (\eta_e\otimes\id)\alpha\big)(\hat W^{(2)})\nonumber \\
& = & \hat \Delta_2\big((\eta_e\otimes \id)(z)\big),
\end{eqnarray}
as well as
\begin{eqnarray*}
(\id \otimes \pi_{\Bo_{G_2}}\otimes \id\otimes \id)\Delta(z)
& = & (\id \otimes \pi_{\Bo_{G_2}}\otimes \id\otimes \id)\big(W_\BG^*(1\otimes z)W_\BG\big)\\
& = & (W^{(1)}\otimes 1)^* (1\otimes z)(W^{(1)}\otimes 1),
\end{eqnarray*}
which implies that
\begin{equation}\label{eqt:pi1-coalg-homo}
(\id \otimes \pi_{\Bo_{G_2}}\otimes \id\otimes \pi_{\Bo_{G_2}})\Delta(z)
= \Delta_1\big((\id \otimes \pi_{\Bo_{G_2}})(z)\big).
\end{equation}

On the other hand, for any $g\in G_1$ and $r\in G_2$, we put $\delta_g$, $\lambda_g$, $\delta_r$ and $\lambda_r$ to be the corresponding elements in $C_0(G_1)$, $C^*(G_1)$, $C_0(G_2)$ and $C^*(G_2)$, respectively.
Note that $W^{(1)} = \sum_{g\in G_1} \delta_g\otimes \lambda_{g}$ and $\hat W^{(2)} = \sum_{t\in G_2} \lambda_t\otimes \delta_{t^{-1}}$.
Hence,
$$(W_{13}^{(1)})^* (1\otimes 1\otimes \delta_g\otimes \lambda_r) W_{13}^{(1)}
\ = \ \sum_{k\in G_1} \delta_k\otimes 1 \otimes \delta_{k^{-1}g} \otimes \lambda_r.$$
Moreover, Relation \eqref{eqt:defn-alpha} implies that $\alpha(\delta_r) = \sum_{h\in G_1}\delta_{h^{-1}} \otimes \delta_{\alpha_h(r)}$, and so
$$(\id\otimes\alpha)\big(\hat W^{(2)}\big)
\ = \ \sum_{h\in G_1}\sum_{t\in G_2} \lambda_{t^{-1}} \otimes \delta_{h^{-1}} \otimes \delta_{\alpha_h(t)}.$$
Consequently,
\begin{eqnarray}\label{eqt:coprod-bicr-prod}
\Delta(\delta_g\otimes \lambda_r)
& = & \sum_{s,t\in G_2}\sum_{f,h,k\in G_1} \delta_k \otimes \lambda_{ts^{-1}}\otimes \delta_h\delta_{k^{-1}g}\delta_f\otimes \delta_{\alpha_{h^{-1}}(t)}\lambda_r \delta_{\alpha_{f^{-1}}(s)}\nonumber \\
& = & \sum_{s,t\in G_2}\sum_{k\in G_1} \delta_k \otimes \lambda_{ts^{-1}}\otimes \delta_{k^{-1}g} \otimes \delta_{\alpha_{g^{-1}k}(t)}\delta_{r \alpha_{g^{-1}k}(s)}\lambda_r \nonumber \\
& = & \sum_{t\in G_2}\sum_{k\in G_1} \delta_k \otimes \lambda_{\alpha_{k^{-1}g}(r)}\otimes \delta_{k^{-1}g} \otimes \delta_{\alpha_{g^{-1}k}(t)}\lambda_r \nonumber \\
& = & \sum_{h\in G_1} \delta_{h^{-1}} \otimes \lambda_{\alpha_{hg}(r)}\otimes \delta_{hg} \otimes \lambda_r.
\end{eqnarray}

Now, supoose that $U\in M(\CK(\KH)\otimes C_0(\BG))$ is a unitary corepresentation of $C_0(\BG)$.
Then
$$U = \sum_{r\in G_2} \sum_{g\in G_1} U_{g,r}\otimes \delta_g\otimes \lambda_r$$
for some $U_{g,r}\in \CL(\KH)$ (here we use $\sum_{g\in G_1} U_{g,r}\otimes \delta_g$ to denote the map $g\mapsto U_{g,r}$).
Relations \eqref{eqt:defn-corep} and \eqref{eqt:coprod-bicr-prod} produce
\begin{eqnarray*}
\lefteqn{\sum_{k,l\in G_1}\sum_{s,t\in G_2} U_{k,s}U_{l,t} \otimes \delta_k\otimes \lambda_s \otimes \delta_l \otimes \lambda_t} \qquad \qquad \qquad \qquad \qquad\\
& = & \sum_{g,h\in G_1}\sum_{r\in G_2} U_{g,r} \otimes \delta_{h^{-1}} \otimes \lambda_{\alpha_{hg}(r)}\otimes \delta_{hg} \otimes \lambda_r\\
& = & \sum_{k,l\in G_1}\sum_{t\in G_2} U_{kl,t} \otimes \delta_{k} \otimes \lambda_{\alpha_{l}(t)}\otimes \delta_l \otimes \lambda_t.
\end{eqnarray*}
This tells us that
\begin{equation}\label{eqt:bicr-unit-corep}
U_{k,s}U_{l,t} = \begin{cases}
U_{kl,t} & \text{ if } s=\alpha_l(t)\\
0        & \text{ otherwise.}
\end{cases}
\end{equation}

Now, we set $U^{(1)}:= (\id\otimes \id \otimes \pi_{\Bo_{G_2}})(U)$ and $U^{(2)}:= (\id\otimes \eta_e\otimes \id)(U)$.
By Relations \eqref{eqt:eta-coalg-homo} and \eqref{eqt:pi1-coalg-homo}, we know that $U^{(1)}$ and $U^{(2)}$ are unitary corepresentations of $C_0(G_1)$ and $C^*(G_2)$, respectively.
They induces, respectively, a unitary representation $\mu_U:G_1\to \CL(\KH)$ and a $^*$-representation $\Psi_U: C(G_2)\to \CL(\KH)$.
Clearly, for any $h\in G_1$ and $t\in G_s$, one has
$$\mu_U(h) = \sum_{r\in G_2} U_{h,r}
\qquad \text{and} \qquad \Psi_U(\delta_t) = U_{e,t}.$$
It is now easy to verify, using Relation \eqref{eqt:bicr-unit-corep}, that
$$\Psi_U(\delta_{\alpha_h(t)})\mu_U(h)  =U_{h,t}= \mu_U(h) \Psi_U(\delta_t)$$
and $(\Psi_U, \mu_U)$ is covariant for the action $\alpha$ and hence induces a $^*$-representation $\Psi_U\times \mu_U: C(G_2)\rtimes_\alpha G_1\to \CL(\KH)$.

On the other hand, as $L^\infty(\BG) = \ell^\infty(G_1)\rtimes_\beta G_2 = \ell^\infty(G_1) \bar{\otimes} L(G_2)$ (see \cite[Definition 3.4.2]{Vaes01} or \cite[Proposition 1.1]{BS98}), its predual
$L^1(\BG)$ is generated by  $\{\bar\lambda_g\otimes \bar\delta_s: g\in G_1; s\in G_2\}$, where $\bar\lambda_g\in \ell^1(G_1)$ and $\bar\delta_s\in A(G_2)$ are the images of $g$ and $s$, respectively.
Thus, 
$$\pi_U(\bar\lambda_g\otimes \bar\delta_s)= U_{g,s}= \mu_U(g) \Psi_U(\delta_s) = (\Psi_U\times \mu_U)(\ti \lambda_g \ti \delta_s),$$ 
where $\ti \delta_s$ and $\ti \lambda_g$ are the images of $\delta_s$ and $\lambda
_g$, respectively, in the full crossed product $C(G_2)\rtimes_\alpha G_1$.

The above shows that $C_0^\ru(\hat \BG)$ is a quotient $C^*$-algebra of $C(G_2)\rtimes_\alpha G_1$.
Since $C(G_2)$ has strong property $T$ (see Example 5.1 and Lemma 4.2 of \cite{LN}), we know, through Theorem 4.6 and Lemma 4.1 of \cite{LN}, that the unital $C^*$-algebra $C_0^\ru(\hat  \BG)$ has strong property $T$, and hence has property $T$.
Now, \cite[Theorem 5.2]{KS} concludes that $\BG$ has property $T$.
\end{prf}

\bigskip

\begin{ack}
The authors would like to thank Prof.\ A.\ Van Daele for some discussion on this topic and for some comments on an earlier version of this work.
\end{ack}

\end{document}